\documentstyle[12pt,twoside]{article}

\font\teneufm=eufm10 scaled \magstep1
\font\seveneufm=eufm7 scaled \magstep1
\font\fiveeufm=eufm5  scaled \magstep1
\newfam\eufmfam
\textfont\eufmfam=\teneufm
\scriptfont\eufmfam=\seveneufm
\scriptscriptfont\eufmfam=\fiveeufm

\newfam\msbfam
\font\tenmsb=msbm10 scaled \magstep1  \textfont\msbfam=\tenmsb
\font\sevenmsb=msbm7 scaled \magstep1 \scriptfont\msbfam=\sevenmsb
\font\fivemsb=msbm5 scaled \magstep1  \scriptscriptfont\msbfam=\fivemsb
\def\Bbb{\fam\msbfam \tenmsb}

\def\RR{{\Bbb R}}
\def\CC{{\Bbb C}}

\def\TT{{\Bbb T}}

\def\ra{\rightarrow}

 \def\HollowBoxx #1#2#3{{\dimen0=#1 \advance\dimen0 by -#2       
       \dimen1=#1 \advance\dimen1 by #3                       
        \vrule height 0pt depth #3 width #2                   
       \hskip -#3
       \vrule height #1 depth #3 width #3}}                   
 \def\LeftContraction{\mathord{\kern1.45pt \HollowBoxx{6pt}{3.5pt}{.4pt}}\,}

 \def\HollowBox #1#2#3{{\dimen0=#1 \advance\dimen0 by -#3       
       \dimen1=#1 \advance\dimen1 by #3                       
        \vrule height #1 depth #3 width #3                    
        \vrule height 0pt depth #3 width #2                   
        \hskip -#3}}                                             
 \def\RightContraction{\mathord{\, \HollowBox{6pt}{3.1pt}{.4pt}} \kern1.6pt}             

\def\qed{{\hfill $\Box$}}
\newtheorem{theorem}{THEOREM}[section]
\newtheorem{lemma}[theorem]{Lemma}

\newtheorem{proposition}[theorem]{Proposition}

\begin{document}
\begin{center}
{\Large \bf Characterization of ${\bf C}^n$
\medskip\\
by its Automorphism Group}\footnote{{\bf Mathematics 
    Subject Classification:} 32Q57, 32M05, 32M17}\footnote{{\bf Keywords 
   and Phrases:} complex manifolds, automorphism groups, equivalence problem.}  
\bigskip\\
\normalsize A. V. Isaev
\end{center}

\begin{quotation} \sl We show that if the group of holomorphic automorphisms 
of a connected Stein manifold $M$ is isomorphic to that of $\CC^n$ as a 
topological group equipped with the compact-open topology, 
then $M$ is biholomorphically equivalent to $\CC^n$.  
\end{quotation} 

\pagestyle{myheadings}
\markboth{A. V. Isaev}
{Characterization of $\CC^n$}

\setcounter{section}{0}

\section{Introduction}
\setcounter{equation}{0}

Let $M$ be a connected complex manifold of dimension $n$ and let
$\hbox{Aut}(M)$ denote the group of holomorphic automorphisms of
$M$. The group $\hbox{Aut}(M)$ is a topological group equipped with
the natural compact-open topology. We are interested in the problem of
characterizing $M$ by $\hbox{Aut}(M)$. This problem becomes
particularly intriguing when $\hbox{Aut}(M)$ is
infinite-dimensional. 

Let, for example, $M=\CC^n$ and suppose that
$M'$ is such that $\hbox{Aut}(M')$ is isomorphic as a topological
group to $\hbox{Aut}(\CC^n)$; is it then true that $M'$ is
biholomorphically equivalent to $\CC^n$?

In \cite{IK} we gave a positive answer to the above question, and
the proof there followed from a general classification of all connected
$n$-dimensional complex manifolds that admit effective actions of the unitary
group $U_n$ by holomorphic transformations. In this paper we give a
simpler proof in the case of Stein manifolds. This proof does not
require considering the whole group $U_n$, but relies only on
linearization of the induced action of the torus $\TT^n\subset U_n$
\cite{BBD}. 

In this paper we prove the following theorem.

\begin{theorem}\label{main} \sl Let $M$ be a connected Stein manifold
  of dimension $n$ and
  suppose that $\hbox{Aut}(M)$ and $\hbox{Aut}(\CC^n)$ are isomorphic
  as topological groups (both groups are considered with the
  compact-open topology). Then $M$ is biholomorphically equivalent to
  $\CC^n$.
\end{theorem}

\section{Proof of Theorem} 
\setcounter{equation}{0}

The theorem is obvious for $n=1$, so we assume that $n\ge 2$. Let $\Phi: \hbox{Aut}(\CC^n)\ra \hbox{Aut}(M)$ be an isomorphism. The
group $\hbox{Aut}(\CC^n)$ contains the subgroup $\CC^{*n}$, i.e.,
all transformations of the form
\begin{equation}
(z_1,\dots,z_n)\mapsto (\lambda_1 z_1,\dots, \lambda_n z_n),\label{dil}
\end{equation}
with $\lambda:=(\lambda_1,\dots,\lambda_n)\in\CC^{*n}$ and $z:=(z_1,\dots,z_n)\in\CC^n$. Therefore $\CC^{*n}$ acts on
$M$ with the action mapping $F:\CC^{*n}\times M\ra M$ defined as follows:
$$
F(\lambda,p):=\Phi(\lambda)(p),
$$
for $\lambda\in\CC^{*n}$, $p\in M$. This action is clearly effective on $M$. Since $F$ is continuous in $(\lambda,p)$ and
holomorphic in $p$, it is in fact real-analytic in
$(\lambda,p)$ \cite{BM}.

We will now prove the following proposition which is of independent
interest and holds for manifolds more general than Stein manifolds.

\begin{proposition}\label{indep}\sl Let $M$ be a connected manifold
  of complex dimension $n$ whose envelope of holomorphy is smooth, and
  suppose that $M$ admits an effective action of $\CC^{*n}$ by
  holomorphic transformations. Then $M$ is
  biholomorphically equivalent to either $\CC^n$, or
  $\CC^n\setminus\{0\}$, or $\CC^n$ without
  some coordinate hyperplanes:
\begin{equation}
\CC^n\setminus\cup_{k=1}^r\{z_{i_k}=0\},\qquad r\ge 1. \label{fin}
\end{equation}
\end{proposition}

\noindent {\bf Proof:} Let as above $F:\CC^{*n}\times M\ra M$ denote the
action of $\CC^{*n}$ on $M$. Consider the restriction of the action to
the torus
$\TT^n:=\{(\lambda_1,\dots,\lambda_n)\in\CC^{*n}:|\lambda_j|=1\,\,\hbox{for
  all $j$}\}$. It follows from \cite{BBD} that 
there exists a holomorphic embedding $\alpha: M\ra \CC^n$ such that
$D:=\alpha(M)$ is a Reinhardt domain and the induced action 
$G(\lambda,z)=(G_1(\lambda,z),\dots,G_n(\lambda,z)):=\alpha(F(\lambda,\alpha^{-1}(z)))$ of $\TT^n$ on $D$ has the form
$$
G_j(\lambda,z)=e^{i(a_{j1}\lambda_1+\dots+a_{jn}\lambda_n)}z_j,
$$
where $\lambda\in\TT^n$, $z\in D$ and
$a_{jk}$ are fixed integers such that $\det(a_{jk})=\pm 1$.

Let $F_{\alpha}:\CC^{*n}\times D\ra D$ be the induced action of $\CC^{*n}$
on $D$:
$$
F_{\alpha}(\lambda,z):=\alpha(F(\lambda,\alpha^{-1}(z))),
$$
for $\lambda\in\CC^{*n}$, $z\in D$. Denote by
$\Phi_{\alpha}:\CC^{*n}\ra\hbox{Aut}(D)$ the corresponding homomorphism:
$$
\Phi_{\alpha}(\lambda)(z)=F_{\alpha}(\lambda,z),
$$
for $\lambda\in\CC^{*n}$, $z\in D$. Since $G$ and $F_{\alpha}$ coincide on
$\TT^n$, it follows that $\Phi_{\alpha}(\TT^n)$ consists precisely of all mappings
of the form (\ref{dil}) with $|\lambda_1|=\dots=|\lambda_n|=1$, i.e., $\Phi_{\alpha}(\TT^n)=\TT^n$. Let
$C(\TT^n)$ denote the centralizer of $\TT^n$ in $\hbox{Aut}(D)$, i.e.,
$$
C(\TT^n):=\{f\in\hbox{Aut}(D):f\circ t=t\circ f\quad\hbox{for all
  $t\in\TT^n$}\}.
$$
Since $\CC^{*n}$ is commutative, we have $\Phi_{\alpha}(\CC^{*n})\subset
C(\TT^n)$.

We now need the following lemma.

\begin{lemma}\label{cent}\sl Every element of $C(\TT^n)$ has the the
  form (\ref{dil}), i.e., $C(\TT^n)\subset\CC^{*n}$.
\end{lemma}

\noindent{\bf Proof:} Let $f=(f_1,\dots,f_n)\in C(\TT^n)$. Then
we have
\begin{equation}
f_j(e^{i\theta_1}z_1,\dots,e^{i\theta_n}z_n)=e^{i\theta_j}f_j(z_1,\dots,z_n),
\qquad j=1,\dots,n,\label{theta}
\end{equation}
for all $\theta_1,\dots,\theta_n\in\RR$ and $z\in D$. In particular, for every
fixed $j$ we have
$$
f_j(z_1,\dots,z_{k-1},e^{i\theta}z_k,z_{k+1},\dots,z_n)=f_j(z_1,\dots,z_n),
$$
for all $\theta\in\RR$, $k\ne j$, which implies that $f_j$ depends
only on $z_j$ (we will write $f_j=f_j(z_j))$. Then (\ref{theta}) gives
\begin{equation}
f_j(e^{i\theta}z_j)=e^{i\theta}f_j(z_j), \label{j}
\end{equation}
for all $\theta\in\RR$. Differentiating (\ref{j}) with respect to
$z_j$ we get
$$
f_j'(e^{i\theta}z_j)=f_j'(z_j),
$$
which implies that $f'_j\equiv \hbox{const}$ and thus $f_j=\lambda_j
z_j$, $\lambda_j\in\CC^*$. 

The lemma is proved.\qed
\smallskip\\

Lemma \ref{cent} gives that $\Phi_{\alpha}$ is a continuous
homomorphism from $\CC^{*n}$ into itself, and thus is a Lie group homomorphism (see, e.g, \cite{W}). Further, since $\Phi_{\alpha}$ is
injective, it is
in fact an automorphism of $\CC^{*n}$. In particular, $D$ is invariant under all
mappings of the form (\ref{dil}), and thus $D$ is either $\CC^n$ or,
$\CC^n\setminus\{0\}$, or a
domain of the form (\ref{fin}).

The proposition is proved.\qed
\smallskip\\

We will now show that the automorphism groups of $\CC^n$ and any
domain of the form (\ref{fin}) are not isomorphic as topological
groups. This is a consequence
of the following observation.

\begin{lemma}\label{connec} \sl For $n\ge 2$ we have
\smallskip\\

\noindent (i) $\hbox{Aut}(\CC^n)$ is connected;
\smallskip\\

\noindent (ii) If $D$ is a domain of the form (\ref{fin}), then
$\hbox{Aut}(D)$ is disconnected.

\end{lemma}

\noindent {\bf Proof:} Following \cite{AL}, we consider special
automorphisms of $\CC^n$ called overshears:
\begin{equation}
(z_1,\dots,z_n)\mapsto
(z_1,\dots,z_{n-1},f(z_1,\dots,z_{n-1})+h(z_1,\dots,z_{n-1})z_n),\label{shear}
\end{equation}
where $f,h$ are entire functions on $\CC^{n-1}$ and $h$ is nowhere
zero. Let $G(\CC^n)$ denote the subgroup of $\hbox{Aut}(\CC^n)$
generated by overshears (\ref{shear}) and permutations of
coordinates. It is proved in \cite{AL} that $G(\CC^n)$ is dense in
$\hbox{Aut}(\CC^n)$. We will show that every element of $G(\CC^n)$ can
be joined with the identity by a continuous path in
$\hbox{Aut}(\CC^n)$. For a mapping of the form (\ref{shear}) we choose
a path $\gamma(t)$ as follows:
\begin{eqnarray*}
&{}&\gamma(t)(z_1,\dots,z_n):=(z_1,\dots,z_{n-1},(1-t)f(z_1,\dots,z_{n-1})+h^{1-t}(z_1,\dots,z_{n-1})z_n),\\
&{}&0\le t\le 1.
\end{eqnarray*}
Further, for the permutation of $z_j$ and $z_k$ we choose (assuming $j<k$)
\begin{eqnarray*}
&{}&\gamma(t)(z_1,\dots,z_n):=\\
&{}&(z_1,\dots,z_{j-1},(1-t)z_k+tz_j,z_{j+1},\dots,z_{k-1},\\
&{}& tz_k+[(1-t)+if(t)]z_j,z_{k+1},\dots,z_n), \qquad 0\le t\le 1,
\end{eqnarray*}
where $f$ is a real-valued continuous function on $[0,1]$ such that $f(0)=f(1)=0$ and
$f(1/2)\ne 0$.

Therefore, $G(\CC^n)$ lies in the identity component of
$\hbox{Aut}(\CC^n)$ and hence $\hbox{Aut}(\CC^n)$ does not in fact
have any other connected components. Thus, $\hbox{Aut}(\CC^n)$ is
connected, and (i) is proved.

Let $D$ be a domain of the form (\ref{fin}). Choose $1\le s\le r$ and
$p\in D$. Let $L_p$ be the complex affine line in $\CC^n$ orthogonal to the hyperplane $\{z_s=0\}$ and passing through
$p$.  Denote by $q$ the point of intersection of $L_p$ and
$\{z_s=0\}$. Next, we choose a closed curve
$\Gamma$ in $L_p\cap D$ around $q$ and define a subset of $\hbox{Aut}(D)$ as:
$$
A:=\left\{(f_1,\dots,f_n)\in\hbox{Aut}(D):\frac{1}{2\pi i} \int_{\Gamma}\frac{{\hat
    f}'_s}{{\hat f}_s}\,dz_s<0\right\},
$$
where ${\hat f}_s:=f_s|_{L_p}$ and ${\hat f}'_s$ denotes the
derivative of ${\hat f}_s$ with respect to $z_s$. 
The subset $A$ is clearly open in $\hbox{Aut}(D)$. It is also closed
in $\hbox{Aut}(D)$
since the integrals in its definition are integers. It is
non-empty since it contains the automorphism
$$
(z_1,\dots,z_n)\mapsto
(z_1,\dots,z_{s-1},\frac{1}{z_s},z_{s+1},\dots,z_n).
$$
And obviously $A\ne \hbox{Aut}(D)$ since the identity mapping is not
in $A$. Therefore, $\hbox{Aut}(D)$ is disconnected, and (ii) is
proved.\qed
\smallskip\\

It now follows from Lemma \ref{connec} that $\hbox{Aut}(\CC^n)$ and
$\hbox{Aut}(D)$ are not isomorphic as topological groups, and hence
$M$ is not equivalent to $D$.  Further, since
$M$ is Stein and $\CC^n\setminus\{0\}$ is not,
Proposition \ref{indep} gives that  $M$ is
biholomorphically equivalent to $\CC^n$, and the theorem is proved. \qed

{\obeylines
Centre for Mathematics and Its Applications 
The Australian National University 
Canberra, ACT 0200
AUSTRALIA 
E-mail address: Alexander.Isaev@anu.edu.au}

\end{document}